\documentclass[a4paper]{article}
\usepackage{amsmath,amssymb,amsthm}
\usepackage{xcolor}
\usepackage{lmodern}
\IfFileExists{upquote.sty}{\usepackage{upquote}}{}
\IfFileExists{microtype.sty}{
  \usepackage[]{microtype}
  \UseMicrotypeSet[protrusion]{basicmath}
}{}
\IfFileExists{xurl.sty}{\usepackage{xurl}}{}
\usepackage{bookmark}
\hypersetup{hidelinks}

\usepackage{amsrefs}

\newtheorem{theorem}{Theorem}[section]
\newtheorem{lemma}[theorem]{Lemma}
\newtheorem{proposition}[theorem]{Proposition}
\newtheorem{corollary}[theorem]{Corollary}
\theoremstyle{definition}
\newtheorem{definition}[theorem]{Definition}

\newtheorem{problem}[theorem]{Problem}

\newcommand{\F}{\mathbb{F}}
\newcommand{\N}{\mathbb{N}}
\newcommand{\Q}{\mathbb{Q}}
\newcommand{\Z}{\mathbb{Z}}

\newcommand{\rad}{\operatorname{rad}}
\newcommand{\divides}{\mathrel{\mid}}
\newcommand{\notdivides}{\mathrel{\nmid}}
\newcommand{\setbuilder}[2]{\left\{#1\,\middle|\,#2\right\}}

\author{Konrad J. Swanepoel\thanks{Department of Mathematics, London School of Economics and Political Science, London WC2A 2AE, United Kingdom. \texttt{k.swanepoel@lse.ac.uk}}}
\title{Equilateral simplices over general fields}
\date{}

\begin{document}
\maketitle

\begin{abstract}
A classical theorem of Schoenberg (1937) characterizes the dimensions~$d$ for which $\Q^d$ contains the vertex set of a regular simplex with $d+1$ vertices.
We generalize Schoenberg's result to all fields of characteristic not equal to~$2$, using elementary Witt theory and Lagrange's four-square theorem.
Our generalization also sharpens Schoenberg's theorem and complements more recent work of Lorenz and Zimmermann (2026) by characterizing the admissible edge lengths.
In the case of finite characteristic, equilateral simplices with $d+2$ vertices can sometimes exist (F.~Rad\'o 1986), and we characterize this situation as well, again for all fields of characteristic not equal to~$2$.
We also briefly consider the edge lengths of equilateral simplices of $d$ vertices.
\end{abstract}

\noindent\textbf{2020 Mathematics Subject Classification.} 11E04 (primary); 52C10, 11E25 (secondary).

\noindent\textbf{Keywords.} Equilateral simplex, quadratic forms over general fields, Schoenberg's theorem, Witt theory, equilateral sets over finite fields.

\section{Introduction}\label{sec:intro}

\subsection{Background and our results}

An old puzzle in combinatorial geometry, going back at least to \'Edouard Lucas (see Beeson~\cite{Beeson} and Khovanova--Radul~\cite{killer}*{Problem~8}), asks whether the integer lattice $\Z^2$ contains the vertex set of an equilateral triangle.
More generally, one may ask: for which dimensions~$d$ does the rational Euclidean space $\Q^d$ contain an equilateral simplex with $d+1$ vertices?

This question was answered by Schoenberg~\cite{Schoenberg1937} and, independently in the even case, by Hadwiger~\cite{Hadwiger1937}, both in 1937.
Schoenberg's characterization states: an equilateral simplex with $d+1$ vertices exists in $\Q^d$ if and only if one of the following holds:
\begin{itemize}
\item $d$ is even and $d+1$ is a perfect square,
\item $d \equiv 3 \pmod{4}$,
\item $d \equiv 1 \pmod{4}$ and $d+1$ is a sum of two integer squares.
\end{itemize}
The first few such dimensions are \[d = 1, 3, 7, 8, 9, 11, 15, 17, 19, 23, 24, 25, 27, 31, 33, 35, 39\ldots,\] recorded as sequence A096315 in the Online Encyclopedia of Integer Sequences~\cite{OEIS_A096315}.
Schoenberg's proof uses the Hasse--Minkowski theory of rational quadratic forms, which classifies rational quadratic forms in terms of local invariants at every prime and at infinity.
The result was rediscovered by Pelling~\cites{Pelling1977, Pelling1980}, who also used the Hasse--Minkowski principle, though to a lesser extent.
Macdonald \cite{Macdonald1987} reformulated Schoenberg's condition by showing it is equivalent to $d+1$ being a sum of $1$, $2$, $4$, or~$8$ odd squares, and gave a proof using quaternionic and octonionic algebras.

In this paper, we generalize Schoenberg's theorem to all fields $F$ of characteristic not equal to~$2$, while simultaneously simplifying the proof.
We define an \emph{equilateral set} in $F^d$ to be a subset $S$ for which there exists a non-zero ``edge length'' $\ell \in F$ (i.e., squared distance) such that $\sum_{i=1}^d(a_i-b_i)^2=\ell$ for all distinct $a,b\in S$.
In accordance with our notation later on, we will denote this space by $(F^d,\sum_{i=1}^d x_i^2)$.
For fields of characteristic~$0$, our generalization is almost identical to Schoenberg's result.
\begin{theorem}\label{cor:schoenberg0}
Let $F$ be a field of characteristic~$0$.
There exists an equilateral set of size $d+1$ in $(F^d, \sum_{i=1}^d x_i^2)$ if and only if:
$d$ is even and $d+1$ is a square in~$F$, or $d \equiv 3 \pmod{4}$, or $d \equiv 1 \pmod{4}$ and $d+1$ is a sum of two squares in~$F$.
\end{theorem}
Specializing to $F = \Q$ and using the classical facts that an integer is a square in~$\Q$ if and only if it is a perfect square, and a positive integer is a sum of two rational squares if and only if it is a sum of two integer squares,
we recover Schoenberg's original result.

Our generalization to fields of finite characteristic is very similar.
Note that when the characteristic of $F$ is $2$, then an equilateral set has size at most~$2$.
We therefore exclude fields of characteristic $2$ from consideration.
This is a consequence of using the sum of squares quadratic form $\sum x_i^2$.
(Lorenz and Zimmermann \cite{LZ2023} consider more general quadratic forms, including over fields of characteristic~$2$.)

\begin{theorem}\label{cor:schoenberg_p}
Let $F$ be a field of characteristic $p$, where $p$ is an odd prime.
There exists an equilateral set of size $d+1$ in $(F^d, \sum_{i=1}^d x_i^2)$ if and only if
$d$ is even and $d+1$ is a square in~$F$ (possibly zero), or $d$ is odd.
Such a set is affinely independent if and only if $p\notdivides d + 1$.
\end{theorem}
When $F$ has finite characteristic, the conditions in Theorem~\ref{cor:schoenberg_p} could have been formulated as in Theorem~\ref{cor:schoenberg0} since $d+1$ is an element of the prime subfield of $F$ where any element is a sum of two squares.
The possibility in Theorem~\ref{cor:schoenberg_p} that an equilateral set of size $d+1$ can lie in a hyperplane was first observed by F.~Rad\'o~\cite{Rado1986}, who also discovered the remarkable fact that in finite characteristic, unit ($\ell=1$) equilateral simplices with $d+2$ vertices can exist in $d$-dimensional space when $p$ divides $d+2$ \cite{Rado1986}*{Propositions~2 and~3}.
This fact was rediscovered many times: see Elsholtz and Klotz~\cite{ElsholtzKlotz}, Krebs \cite{Krebs}*{Section~5}, Lorenz and Zimmermann \cite{LZ2023}*{Lemma~4.9} and Solymosi and Yip \cite{Solymosi2025}.

We give a characterization in any field of finite characteristic not equal to $2$.
\begin{theorem}\label{cor:d+2}
Let $F$ be a field of characteristic $p$, where $p$ is an odd prime.
There exists an equilateral set of size $d+2$ in $(F^d, \sum_{i=1}^d x_i^2)$ if and only if $p\divides d + 2$ and either $d$ is even and $-1$ is a square in~$F$, or $d$ is odd.
\end{theorem}
In particular, for finite fields $\F_{p^k}$ we have the following characterization, since $-1$ is a square in $\F_{p^k}$ if and only if $p \equiv 1 \pmod{4}$ or $k$ is even. This finite-field characterization also appears in Solymosi and Yip \cite{Solymosi2025}*{Proposition~4.4}.

\begin{corollary}\label{cor:schoenberg_finite}
Let $q = p^k$ with $p$ an odd prime.
There exists an equilateral set of size $d+2$ in $(\F_q^d, \sum_{i=1}^d x_i^2)$ if and only if $p \mid d+2$ and either $d$ is odd, or $d$ is even and $k$ is even, or $d$ is even and $p \equiv 1 \pmod{4}$.
\end{corollary}

Since the Hasse--Minkowski principle applies only over global fields such as~$\Q$, and not over arbitrary fields of characteristic~$0$, the methods of Schoenberg~\cite{Schoenberg1937} and Pelling~\cite{Pelling1977} do not extend to this generality. Our proofs use elementary Witt theory of quadratic forms over fields of characteristic not equal to $2$, specifically, diagonalizability of quadratic forms and Witt's cancellation and extension theorems, together with the two-square and four-square composition identities and, in characteristic~$0$, Lagrange's four-square theorem.

To prove Theorems~\ref{cor:schoenberg0}, \ref{cor:schoenberg_p}, and \ref{cor:d+2}, we first sharpen these results by characterizing, for each dimension~$d$, which non-zero $\ell$ allow an $\ell$-equilateral set of size $d+1$ or $d+2$ to exist, in Theorems~\ref{thm:schoenberg0}, \ref{thm:schoenbergp}, and \ref{cor:rado}.

\subsection{Related work}\label{section:related_work}
Chilakamarri~\cite{Chilakamarri1988} studied unit-distance graphs in $\Q^d$ and gave Diophantine conditions for the existence of unit-length equilateral sets.
Their approach uses a result of Hall and Ryser~\cite{HallRyser} that is a special case of the Witt extension theorem.
Elsholtz and Klotz~\cite{ElsholtzKlotz} extended these results to arbitrary fields of characteristic not~$2$, giving conditions for equilateral sets of various sizes, though their nonexistence criterion for unit-equilateral sets of size~$d$ in odd dimension assumes that the quadratic form $\sum x_i^2$ is anisotropic, i.e., $\sum_i x_i^2=0$ only if $x_1=x_2=\dots=x_d=0$.
In characteristic~$0$, Lorenz and Zimmermann~\cite{LZ2023}*{Theorem~6.5} obtained this characterization without the anisotropy assumption, using the full Witt theory of quadratic forms; they also considered arbitrary non-degenerate quadratic forms.
Krebs~\cite{Krebs} determined the maximum unit-length equilateral set sizes over finite fields.
Macdonald \cite{Macdonald1987} and Maehara and Martini~\cite{MaeharaMartini2018} gave elementary constructions of equilateral rational simplices.
See also \cites{Sokolov,Ionascu2007,Ionascu2008a,Ionascu2009,Ionascu2011a,Ionascu2015,Noble2022} for related work on equilateral simplices.

Kominers~\cite{Kominers2026}*{Theorem~1} characterizes the positive rational squared edge lengths of regular simplices in $\Q^n$ in every dimension and codimension, using Hasse--Minkowski theory.
His full-dimensional case agrees with the special case of Theorem~\ref{thm:schoenberg0} below where $F=\Q$, and his codimension-one results settle the rational case of our concluding Problem~\ref{problem}.
The present paper treats arbitrary fields of characteristic not equal to~$2$, including affinely dependent equilateral sets in finite characteristic, using elementary Witt theory.

\subsection{Organisation}
In Section~\ref{sec:prelim} we recall the necessary background on quadratic forms.
Section~\ref{sec:equilateral_form} introduces the equilateral quadratic form~$q_d$ and relates it to equilateral sets via a lemma of Schoenberg.
Section~\ref{subsec:singular} treats the singular case that arises in finite characteristic, where affinely dependent equilateral sets can exist.
Section~\ref{sec:structure} develops the structural properties of~$q_d$ needed for the main theorems, treating characteristic~$0$ and finite characteristic separately.
Section~\ref{sec:schoenberg} contains the proofs of Theorems~\ref{cor:schoenberg0}, \ref{cor:schoenberg_p}, and~\ref{cor:d+2}.
Finally, in Section~\ref{sec:size_d} we discuss our attempts to characterize equilateral simplices of given edge length with $d$ vertices in a similar way as in Theorems~\ref{thm:schoenberg0} and \ref{thm:schoenbergp}.

\section{Preliminaries on quadratic forms}\label{sec:prelim}

We recall the necessary background on quadratic forms; all of the following can be found in Lam~\cite{Lam}*{Chapter~I}.

Throughout this paper, $F$ denotes a field of characteristic $\chi(F) \neq 2$.
We write $F^*$ for the multiplicative group $F \setminus \{0\}$, $F^{*2}$ for the subgroup of squares $\setbuilder{x^2}{x \in F^*}$, and $F^*/F^{*2}$ for the group of \emph{square classes}.

A \emph{quadratic space} over~$F$ is a pair $(V, B) = (V, q)$ where $V$ is a finite-dimensional vector space over~$F$, $B \colon V \times V \to F$ is a symmetric bilinear form, and $q(x) = B(x,x)$ is the associated quadratic map.
Since $\chi(F) \neq 2$, the bilinear form can be recovered from $q$ by polarization: $B(x,y) = \frac{1}{2}(q(x+y) - q(x) - q(y))$.
With respect to a basis $e_1, \ldots, e_d$ of~$V$, the bilinear form is represented by a symmetric matrix~$M$ via $B(x,y) = x^T M y$, and $q$ becomes a homogeneous polynomial of degree~$2$.

Two quadratic forms $q_1$ and $q_2$ are \emph{equivalent} (or \emph{isometric}), written $q_1 \cong q_2$, if their matrices $M_1$ and $M_2$ are \emph{congruent}, i.e., there exists an invertible matrix~$A$ such that $M_1 = A^T M_2 A$.
A quadratic space is \emph{non-degenerate} if $B(x,y) = 0$ for all~$y$ implies $x = o$, equivalently, if the matrix~$M$ is invertible.
The \emph{determinant} of a non-degenerate form is the square class $\det(M) F^{*2}$, which is well defined up to equivalence as congruent matrices have the same determinants up to a square factor.

If $S$ is a subspace of $(V,B)$, we denote its \emph{orthogonal complement} by $S^\perp=\setbuilder{x\in V}{B(x,s)=0\;\forall s\in S}$, and the \emph{radical} of $V$ by $\rad V:= V^\perp$.
We will use the fact that if $S$ is a non-degenerate subspace of $V$ then $V=S\perp S^\perp$~\cite{Lam}*{Chapter~I, Corollary~2.5}.

We say $a \in F^*$ is \emph{represented} by~$q$ if $a = q(v)$ for some $v \in V$, and write $D(q)$ for the set of non-zero values represented by~$q$.
By the \emph{representation criterion}~\cite{Lam}*{Chapter~I, 2.3}, $a \in D(q)$ if and only if $q \cong \langle a \rangle \perp q'$ for some form~$q'$.

We write $\langle a_1, \ldots, a_d \rangle$ for the \emph{diagonal form} $\sum_{i=1}^d a_i x_i^2$, and $k \cdot \langle a \rangle$ for the diagonal form
\(\sum_{i=1}^k ax_i^2\).
For $a \in F^*$, we write $a \times q$ for the quadratic form $x \mapsto a q(x)$.
The \emph{orthogonal sum} of quadratic spaces $(V_1, q_1)$ and $(V_2, q_2)$ is denoted $q_1 \perp q_2$.
Every quadratic form over~$F$ can be diagonalized~\cite{Lam}*{Chapter~I, Corollary~2.4}.

Two binary non-degenerate forms $\langle a, b \rangle$ and $\langle c, d \rangle$ are isometric if and only if they have the same determinant and there exists $e \in F^*$ represented by both~\cite{Lam}*{Chapter~I, Proposition~5.1}.
We note two consequences.
If $a$ is a non-zero sum of two squares in~$F$, then $\langle a, a \rangle \cong \langle 1, 1 \rangle$.
The \emph{hyperbolic plane} $\langle 1, -1 \rangle$ satisfies $\langle 1, -1 \rangle \cong \langle a, -a \rangle$ for all $a \in F^*$, since $a=(\frac{a+1}{2})^2-(\frac{a-1}{2})^2$.

We say $a \in F^*$ is a \emph{similarity factor} of~$q$ if $q \cong a \times q$.
The set $G(q)$ of similarity factors is a subgroup of~$F^*$ containing~$F^{*2}$.
If $q$ is non-degenerate and its dimension is odd, then $G(q) = F^{*2}$.
The two-square and four-square composition identities, valid over~$F$, show that the forms $2 \cdot \langle 1 \rangle$ and $4 \cdot \langle 1 \rangle$ satisfy $D(q)=G(q)$ and hence are \emph{group forms}: $D(q)$ is a subgroup of~$F^*$.
In particular, $2 \cdot \langle 1 \rangle \cong 2 \cdot \langle a \rangle$ for any non-zero~$a$ that is a sum of two squares in~$F$, and $4 \cdot \langle 1 \rangle \cong 4 \cdot \langle a \rangle$ for any non-zero~$a$ that is a sum of four squares in~$F$.

A non-zero vector $v$ in $(V,q)$ is \emph{isotropic} if $q(v)=0$,
and anisotropic if $q(v)\neq 0$. The space $(V,q)$ is
\emph{isotropic} if it contains a non-zero isotropic vector, and
\emph{anisotropic} if all non-zero vectors are anisotropic. The space
$(V,q)$ is \emph{totally isotropic} if all non-zero vectors in $V$
are isotropic. This is equivalent to $B$ being the zero bilinear form.

We rely on the following two fundamental results.

\begin{theorem}[Witt cancellation, {\cite{Lam}*{Chapter~I, 4.2}}]\label{thm:witt_cancel}
For any quadratic forms $q, q_1, q_2$, if $q \perp q_1 \cong q \perp q_2$, then $q_1 \cong q_2$.
\end{theorem}

\begin{theorem}[Witt extension, {\cite{Lam}*{Chapter~I, Exercise~11}}]\label{thm:witt_ext}
Let $V$ be a non-degenerate quadratic space with subspaces $S_1$ and $S_2$.
Then any isometry $S_1 \to S_2$ can be extended to an isometry $V \to V$.
\end{theorem}

\section{The equilateral quadratic form}\label{sec:equilateral_form}

\begin{definition}
Let $\ell \in F^*$.
A subset~$S$ of a quadratic space $(V, q)$ is \emph{$\ell$-equilateral} if $q(a - b) = \ell$ for all distinct $a, b \in S$.
\end{definition}

Intuitively, an $\ell$-equilateral set is the vertex set of a regular simplex with squared side length~$\ell$.
We work with the squared length~$\ell$ rather than $\sqrt{\ell}$ to avoid issues with the existence of square roots.
We exclude $\ell = 0$, since the affine span of a set $S$ satisfying $q(a-b)=0$ for all  $a,b\in S$ is a translate of a totally isotropic subspace.
(This can be seen in the same way as the proof of the next lemma.)

For each $d$ we define the quadratic form
\[
q_d(X_1, \ldots, X_d) = \sum_{i=1}^d \sum_{j=i}^d X_i X_j,
\]
whose matrix $M_d$ has diagonal entries~$1$ and off-diagonal entries~$1/2$.
It is easily seen that $\det(M_d)=(d+1)/2^d$, so $q_d$ is non-degenerate if and only if $d + 1 \neq 0$ in~$F$.
In this case,
\[
\det(q_d) = \begin{cases} d+1 & \text{if } d \text{ is even,} \\ (d+1)/2 & \text{if } d \text{ is odd,} \end{cases}
\quad\text{(modulo $F^{*2}$).}
\]

The next lemma is due to Schoenberg (see also \cite{LZ2023}*{Lemma~4.5}).

\begin{lemma}[Schoenberg~\cite{Schoenberg1937}]\label{lem:gram}
Let $\ell\in F^*$.
Then a $d$-dimensional quadratic space $(V, q)$ contains an affinely independent $\ell$-equilateral set of size $d+1$ if and only if
$q$ is isometric to $\ell \times q_d$.
If $(V,q)$ is non-degenerate, this is possible only if $d+1\neq 0$ in~$F$.
\end{lemma}
\begin{proof}
Let $\{a_0, a_1, \ldots, a_d\}$ be the equilateral set, and set $u_i = a_i - a_0$ for $i = 1, \ldots, d$.
Then $u_1, \ldots, u_d$ form a basis of~$V$, with $B(u_i, u_i) = q(a_i - a_0) = \ell$.
For distinct $i, j$,
\[
\ell = q(a_i - a_j) = q(u_i - u_j) = 2\ell - 2B(u_i, u_j),
\]
hence $B(u_i, u_j) = \ell/2$.
Thus, in this basis, $q$ has matrix $\ell M_d$, so $q\cong\ell\times q_d$.
Conversely, the origin and the standard basis vectors form an affinely independent $\ell$-equilateral set in $(F^d,\ell\times q_d)$, and an isometry transports this set to $(V,q)$.
Finally, if $(V,q)$ is non-degenerate, then $0\neq\det(\ell M_d)=\ell^d(d+1)/2^d$, so $d+1\neq0$ in~$F$.
\end{proof}

The following lemma shows that affinely independent equilateral sets are essentially unique if the space is non-degenerate.

\begin{lemma}\label{lem:unique}
If $(V,q)$ is non-degenerate, then any two affinely independent $\ell$-equilateral sets of the same cardinality can be mapped to each other by an affine isometry $x\mapsto T(x)+b$, where $T$ is an isometry of~$(V,q)$ and $b\in V$.
\end{lemma}
\begin{proof}
Translate a point of each set to the origin. The resulting difference vectors have the same Gram matrix and are linearly independent, so the linear map matching them is an isometry of their spans. Apply the Witt extension theorem (Theorem~\ref{thm:witt_ext}) and undo the translations.
\end{proof}

\section{Affinely dependent equilateral sets}\label{subsec:singular}

In fields of finite characteristic, it is possible for affinely dependent equilateral sets to exist, as discovered by Ferenc Rad\'o~\cite{Rado1986} for finite fields.
See also \cites{Krebs, LZ2023, Solymosi2025}.

\begin{lemma}\label{lem:rado}
Let $\ell\in F^*$.
Let $S$ be an affinely dependent $\ell$-equilateral set of size $k+1$ in $(V, q)$.
Then $k + 1 = 0$ in~$F$, any $k$ points of~$S$ are affinely independent, and each point of~$S$ is the centroid of the remaining~$k$ points.
Conversely, if $k + 1 = 0$ in~$F$, then any affinely independent $\ell$-equilateral set of size~$k$ can be extended to an affinely dependent $\ell$-equilateral set of size $k+1$ by adding its centroid.
\end{lemma}
\begin{proof}
By translation, we may assume $o \in S$.
The Gram matrix of the remaining $k$ vectors is $\ell M_k$.
Since $S$ is affinely dependent, this matrix must be singular, so $(k+1)(\ell/2)^k = 0$ in~$F$, giving $k + 1 = 0$ in~$F$ (since $\ell \neq 0$).
The same argument applied to any $k$ points of~$S$ shows they cannot be affinely dependent, since $k = -1 \neq 0$ in~$F$.

Consider the points $o, e_1, \ldots, e_{k-1}\in S$ and write the remaining point of $S$ as $\sum \lambda_i e_i$.
We compute from $q(\sum \lambda_i e_i) = \ell$ and $q(e_j - \sum \lambda_i e_i) = \ell$ for all~$j$, together with $k + 1 = 0$, that $\lambda_i = -1 = 1/k$ for all~$i$, so $\sum \lambda_i e_i$ is the centroid of $\{o, e_1, \ldots, e_{k-1}\}$.
Since the dependent point is uniquely determined from the others, each point $s\in S$ is the centroid of $S\setminus\{s\}$.

The converse follows by checking that the centroid of an $\ell$-equilateral set of size~$k$ is at squared distance~$\ell$ from each point when $k + 1 = 0$ in~$F$.
\end{proof}

Lemmas~\ref{lem:gram} and \ref{lem:rado} give the following.
\begin{theorem}\label{cor:rado}
Let $\ell\in F^*$.
The following are equivalent:
\begin{itemize}
\item There exists an $\ell$-equilateral set of size $d+2$ in $(F^d, q)$.
\item $d + 2 = 0$ in~$F$ and there exists an affinely independent $\ell$-equilateral set of size $d+1$ in $(F^d,q)$.
\item $d + 2 = 0$ in~$F$ and $q \cong \ell \times q_d$.
\end{itemize}
In this case, any point of the $d+2$ points is the centroid of the other $d+1$.
There does not exist an $\ell$-equilateral set of size $d+3$ in $(F^d, q)$.
\qed
\end{theorem}

\begin{lemma}\label{lem:unique2}
If $(V,q)$ is non-degenerate, then any two affinely dependent $\ell$-equilateral sets of the same cardinality can be mapped to each other by an affine isometry of~$(V,q)$.
\end{lemma}
\begin{proof}
This is immediate from Lemma~\ref{lem:rado} and Lemma~\ref{lem:unique} applied to the two affinely independent sets obtained by removing a point from each of the affinely dependent equilateral sets.
\end{proof}
In Section~\ref{subsec:charp} we characterize the quadratic forms isometric to $q_d$ when the characteristic is finite.

\section{Structural results on the equilateral form}\label{sec:structure}

In this section we develop the key properties of~$q_d$ needed for the main theorems.
We consider the cases of characteristic $0$ and finite characteristic separately.

\subsection{Characteristic zero}

\begin{lemma}\label{lem:decomp}
If $F$ has characteristic $0$, then for any $d \geq 1$,
$q_d \cong q_{d-1} \perp \langle \tfrac{1}{2} d(d+1) \rangle$.
\end{lemma}
\begin{proof}
The case $d=1$ is immediate, with $q_0$ the zero-dimensional form.
For $d\geq 2$, let $H=\setbuilder{(x_1,\ldots,x_d)\in F^d}{x_d=0}$.
The restriction of $q_d$ to $H$ is $q_{d-1}$, which is non-degenerate since $\det(M_{d-1})=d/2^{d-1}\neq 0$.
Thus $F^d=H\perp H^\perp$, and since $q_d$ is non-degenerate and $H^\perp$ is one-dimensional, $q_d\cong q_{d-1}\perp\langle a\rangle$ for some $a\in F^*$.
Comparing determinants as square classes gives
\[
aF^{*2}
=\frac{\det(M_d)}{\det(M_{d-1})}F^{*2}
=\frac{d+1}{2d}F^{*2}
=\frac{d(d+1)}{2}F^{*2},
\]
where the last equality uses the non-zero square $d^2$.
Hence $\langle a\rangle\cong\langle\tfrac{1}{2}d(d+1)\rangle$, as required.
\end{proof}

From Lemma~\ref{lem:decomp} and induction, we obtain the following result of Schoenberg (see also \cite{LZ2023}*{Lemma~6.1}).

\begin{proposition}[Schoenberg~\cite{Schoenberg1937}]\label{prop:diag}
If $F$ has characteristic~$0$, then
\[
q_d \cong \left\langle \frac{2}{2}, \frac{3}{4}, \frac{4}{6}, \frac{5}{8}, \ldots, \frac{d+1}{2d} \right\rangle.
\]
\end{proposition}
Schoenberg~\cite{Schoenberg1937} proved this diagonalization of $q_d$ by observing that the following identity holds:
$q_d(x_1, \ldots, x_d) = \frac{1}{2} \sum_{i=1}^d \frac{i+1}{i} (x_i + \frac{1}{i+1} \sum_{j > i} x_j)^2$.
Lorenz and Zimmermann \cite{LZ2023}*{Lemma~6.1} gave a matrix argument.

In a field of characteristic \(0\), it is known \cite{Lam}*{Chapter II, Exercise 9}, \cite{AMM1936}*{Editorial Note, p.~254} that the quadratic form
\(-\sum_{i,j=1}^{d+2}|i-j|x_ix_j\) is isometric to
\(\langle 1,-1,\frac{2}{2},\frac{3}{4},\frac{4}{6},\dots,\frac{d+1}{2d}\rangle\cong\langle 1,-1\rangle\perp q_d\),
and also isometric to
\((d+1)\cdot\langle 2\rangle\perp\langle -2(d+1)\rangle\) \cite{Lam}*{Chapter II, Exercise 7(c)}, \cite{Lam1999}*{p.~234}. Thus,
\((d+1)\cdot\langle 2\rangle\perp\langle-2(d+1)\rangle\cong\langle1,-1\rangle\perp q_d\).
We present a direct proof of this consequence.

\begin{proposition}\label{prop:hyperbolic}
If $F$ has characteristic $0$, then for any $d \geq 0$,
\[
q_d \perp \langle 1, -1 \rangle \cong (d+1) \cdot \langle 2 \rangle \perp \langle -2(d+1) \rangle.
\]
\end{proposition}
\begin{proof}
For $d = 0$, both forms are the hyperbolic plane $\langle 1, -1 \rangle \cong \langle 2, -2 \rangle$.
For $d > 0$, we proceed by induction.
Assuming $q_{d-1} \perp \langle 1, -1 \rangle \cong d \cdot \langle 2 \rangle \perp \langle -2d \rangle$ and using $q_d \cong q_{d-1} \perp \langle \frac{1}{2} d(d+1) \rangle$ from Lemma~\ref{lem:decomp}, by Witt cancellation we only need to show that $\langle -2d, \frac{1}{2} d(d+1) \rangle \cong \langle 2, -2(d+1) \rangle$.
This follows as the two forms have the same determinant, and $-2d = 2 - 2(d+1)$ is represented by both.
\end{proof}

\begin{corollary}\label{cor:three_copies}
If $F$ has characteristic~$0$, then for any $d \geq 3$,
\[
q_d \cong (d-3) \cdot \langle 2 \rangle \perp 3 \cdot \langle 2(d+1) \rangle.
\]
\end{corollary}
\begin{proof}
Since $D(4 \cdot \langle 1 \rangle)=G(4 \cdot \langle 1 \rangle)$ and, by Lagrange's four-square theorem, any positive integer is a sum of $4$ squares, $\langle 2, 2, 2, 2 \rangle \cong 4 \cdot \langle 2(d+1) \rangle$.
Also, $\langle 1, -1 \rangle \cong \langle 2(d+1), -2(d+1) \rangle$, since both are the hyperbolic plane.
The result follows from Proposition~\ref{prop:hyperbolic} and Witt cancellation.
\end{proof}

The following characterizes when $q_d$ is isometric to the standard form $d \cdot \langle 1 \rangle$.

\begin{theorem}\label{thm:type0}
Let $F$ have characteristic~$0$.
Then $q_d \cong d \cdot \langle 1 \rangle$ if and only if $\det(q_d)$ is a square in~$F$.
\end{theorem}
\begin{proof}
Necessity is clear from $\det(d \cdot \langle 1 \rangle) = 1$.

Next assume $\det(q_d)$ is a square.
If $d$ is odd, then $\det(q_d) = 2(d+1)$ up to a square, so by Proposition~\ref{prop:hyperbolic} and $\langle 2, 2 \rangle \cong \langle 1, 1 \rangle$,
\[
q_d \perp \langle 1, -1 \rangle \cong (d+1) \cdot \langle 2 \rangle \perp \langle -1 \rangle \cong (d+1) \cdot \langle 1 \rangle \perp \langle -1 \rangle.
\]
Witt cancellation gives $q_d \cong d \cdot \langle 1 \rangle$.

If $d$ is even, then $\det(q_d) = d+1$ up to a square, and similarly
\[
q_d \perp \langle 1, -1 \rangle \cong (d+1) \cdot \langle 2 \rangle \perp \langle -2 \rangle \cong d \cdot \langle 1 \rangle \perp \langle 2, -2 \rangle.
\]
Since $\langle 2, -2 \rangle \cong \langle 1, -1 \rangle$, the result again follows from Witt cancellation.
\end{proof}

\subsection{Characteristic \texorpdfstring{$p$}{p}}\label{subsec:charp}

Over finite prime fields, every non-degenerate quadratic form is determined up to isometry by its dimension and determinant, since the group of square classes has order~$2$ and $\langle \nu, \nu \rangle \cong \langle 1, 1 \rangle$ for any non-square~$\nu$, because every element of the prime subfield is a sum of two squares.
This gives the following.

\begin{lemma}\label{lem:fin_qd}
Let $d \in \N$ and let $F$ be a field of finite characteristic $p \neq 2$.
Then
\[
q_d \cong \begin{cases}
(d-1) \cdot \langle 1 \rangle \perp \langle \det(q_d) \rangle & \text{if } d + 1 \neq 0 \text{ in } F, \\
(d-2) \cdot \langle 1 \rangle \perp \langle \det(q_{d-1}), 0 \rangle & \text{if } d + 1 = 0 \text{ in } F.
\end{cases}
\]
\end{lemma}
\begin{proof}
If $d + 1 \neq 0$ in~$F$, then $q_d$ is non-degenerate.
Since the matrix of $q_d$ has entries in the prime subfield, we can
diagonalize over the prime subfield to obtain $q_d \cong \langle a_1, \ldots, a_d \rangle$ where each $a_i \in \{1, \nu\}$ (with $\nu$ non-square modulo~$p$).
Since $\langle \nu, \nu \rangle \cong \langle 1, 1 \rangle$, we can assume that at most one~$a_i$ equals~$\nu$.
Matching the determinant gives the result.

If $d + 1 = 0$ in~$F$, then $q_d$ is degenerate.
Let $v$ be a non-zero element of the radical of $(F^d,q_d)$, i.e., $v\neq o$ and $B(x,v)=0$ for all $x\in F^d$.
Let $S$ be a subspace of codimension $1$ not containing~$v$.
The projection along~$v$ sends any $1$-equilateral set of size $d+1$ to an affinely dependent equilateral set in the $(d-1)$-dimensional space~$S$, which is thus isometric to $q_{d-1}$ by Theorem~\ref{cor:rado}.
Since $p \nmid d$, the form $q_{d-1}$ is non-degenerate, so $q_d \cong q_{d-1} \perp \langle 0 \rangle$.
\end{proof}
As a consequence, we obtain the finite characteristic analogue of Theorem~\ref{thm:type0}.
\begin{theorem}\label{thm:type0p}
Let $F$ have finite characteristic.
Then $q_d \cong d \cdot \langle 1 \rangle$ if and only if $\det(M_d)=(d+1)/2^d$ is a non-zero square in~$F$.
\end{theorem}
\begin{proof}
Necessity is clear.
For sufficiency, $q_d \cong (d-1) \cdot \langle 1 \rangle \perp \langle \det(q_d) \rangle$ by Lemma~\ref{lem:fin_qd}, and if $\det(q_d)$ is a square, then $\langle \det(q_d) \rangle \cong \langle 1 \rangle$.
\end{proof}

\section{Proofs of the main results}\label{sec:schoenberg}

We now state and prove the main results, our generalizations of Schoenberg's theorem to arbitrary fields.
We treat the cases of zero and finite characteristic separately, as the technical details have different characters.

\subsection{Characteristic zero}

\begin{theorem}\label{thm:schoenberg0}
Let $F$ be a field of characteristic~$0$, and $\ell \in F^*$.
There exists an $\ell$-equilateral set of size $d+1$ in $(F^d, \sum_{i=1}^d x_i^2)$ if and only if:
\begin{itemize}
\item when $d$ is even: $d+1$ is a square in~$F$ and $\ell$ is a similarity factor of $d \cdot \langle 1 \rangle$;
\item when $d \equiv 3 \pmod{4}$: $\ell = (d+1)/2$ up to a square in~$F$;
\item when $d \equiv 1 \pmod{4}$: $d+1$ is a sum of two squares in~$F$, and $\ell = (d+1)/2$ up to a square in~$F$.
\end{itemize}
\end{theorem}
\begin{proof}
By Lemma~\ref{lem:gram}, the existence of such a set is equivalent to $\ell \times q_d \cong d \cdot \langle 1 \rangle$.
For $d=1$, this holds if and only if $\ell$ is a square, as required, so we assume that $d\geq 2$.

\smallskip First consider the case where $d$ is even.
If $\ell \times q_d \cong d \cdot \langle 1 \rangle$, comparing determinants gives $\ell^d (d+1) = 1$ up to a square, so $d + 1$ is a square.
By Theorem~\ref{thm:type0}, $q_d \cong d \cdot \langle 1 \rangle$, so $\ell$ is a similarity factor of $d \cdot \langle 1 \rangle$.
The converse follows by reversing the argument.

\smallskip Next, consider the case $d \equiv 3 \pmod{4}$.
Suppose that
$\ell\times q_d\cong d\cdot\langle 1\rangle$.
Comparing determinants, we get that \(\ell=(d+1)/2\) up to a square.

Conversely, suppose $\ell = (d+1)/2$ up to a square.
By Corollary~\ref{cor:three_copies},
$\ell \times q_d \cong (d-3) \cdot \langle d+1 \rangle \perp \langle 1, 1, 1 \rangle$.
Since $d - 3$ is divisible by~$4$ and $4 \cdot \langle 1 \rangle\cong 4\cdot\langle d+1\rangle$, $(d-3) \cdot \langle d+1 \rangle \cong (d-3) \cdot \langle 1 \rangle$.
Thus $\ell \times q_d \cong d \cdot \langle 1 \rangle$.

\smallskip Finally, suppose that $d \equiv 1 \pmod{4}$.
As in the preceding case, if $\ell\times q_d\cong d\cdot\langle 1\rangle$,
then \(\ell=(d+1)/2\) up to a square, giving $\ell \times q_d \cong (d-3) \cdot \langle d+1 \rangle \perp \langle 1, 1, 1 \rangle$.
Since $d - 5$ is divisible by~$4$, we get $(d-5) \cdot \langle d+1 \rangle \cong (d-5) \cdot \langle 1 \rangle$, so $\ell \times q_d \cong (d-2) \cdot \langle 1 \rangle \perp \langle d+1, d+1 \rangle$.
By Witt cancellation, this is $\cong d \cdot \langle 1 \rangle$ only if $\langle d+1, d+1 \rangle \cong \langle 1, 1 \rangle$, which is equivalent to $d+1$ being a sum of two squares in~$F$.
The converse is obtained by reversing the argument.
\end{proof}

\begin{proof}[Proof of Theorem~\ref{cor:schoenberg0}]
When $d$ is even, we use $\ell = 1$ in Theorem~\ref{thm:schoenberg0}, and when $d$ is odd, we use $\ell=(d+1)/2$.
\end{proof}

\subsection{Finite characteristic}
In finite characteristic we need to distinguish between affinely independent and affinely dependent equilateral sets.

\begin{theorem}\label{thm:schoenbergp}
Let $F$ be a field of characteristic $p$, where $p$ is an odd prime, and let $\ell \in F^*$.
\begin{enumerate}
\item\label{item:aff_indep} There exists an affinely independent $\ell$-equilateral set of size $d+1$ in $(F^d, \sum_{i=1}^d x_i^2)$ if and only if:
  \begin{itemize}
  \item when $d$ is even: $d+1$ is a non-zero square in~$F$ and $\ell$ is a similarity factor of $d \cdot \langle 1 \rangle$;
  \item when $d$ is odd: $d + 1 \neq 0$ in~$F$, and $\ell = (d+1)/2$ up to a square in~$F$.
  \end{itemize}
\item\label{item:aff_dep} There exists an affinely dependent $\ell$-equilateral set of size $d+1$ in $(F^d, \sum_{i=1}^d x_i^2)$ if and only if $d + 1 = 0$ in~$F$ and $\ell$ is a similarity factor of
  \begin{itemize}
  \item $d \cdot \langle 1 \rangle$ when $d$ is even,
  \item $(d+1) \cdot \langle 1 \rangle$ when $d$ is odd.
  \end{itemize}
\end{enumerate}
\end{theorem}
\begin{proof}
Statement \ref{item:aff_indep}.\ when $d$ is even is proved exactly as in Theorem~\ref{thm:schoenberg0}, using Theorem~\ref{thm:type0p} instead of Theorem~\ref{thm:type0}.

\smallskip
Consider statement~\ref{item:aff_indep}.\ when $d$ is odd.
If $\ell \times q_d \cong d \cdot \langle 1 \rangle$, then $q_d$ is non-degenerate, so $d + 1 \neq 0$.
Comparing determinants gives $\ell = (d+1)/2$ up to a square.
Conversely, if $d + 1 \neq 0$ and $\ell = (d+1)/2$ up to a square, then by Lemma~\ref{lem:fin_qd},
$q_d \cong (d-1) \cdot \langle 1 \rangle \perp \langle (d+1)/2 \rangle$,
so $\ell \times q_d \cong (d-1) \cdot \langle (d+1)/2 \rangle \perp \langle 1 \rangle$.
Since $\langle (d+1)/2, (d+1)/2 \rangle \cong \langle 1, 1 \rangle$ ($(d+1)/2$ is in the prime subfield), we conclude $\ell \times q_d \cong d \cdot \langle 1 \rangle$.

\smallskip
Next, consider statement~\ref{item:aff_dep}.\ when $d$ is even.
Suppose there is an affinely dependent $\ell$-equilateral set of size $d+1$.
By Lemma~\ref{lem:rado}, $d + 1 = 0$ in~$F$, and by Lemma~\ref{lem:gram} a translate of the equilateral set spans a hyperplane isometric to $\ell \times q_{d-1}$, which is isometric to $(d-2) \cdot \langle \ell \rangle \perp \langle \ell d/2 \rangle$ by Lemma~\ref{lem:fin_qd}.
Since this hyperplane is non-degenerate and embeds in $d \cdot \langle 1 \rangle$, we compute
$d \cdot \langle 1 \rangle \cong (d-2) \cdot \langle \ell \rangle \perp \langle \ell d/2, \ell d/2 \rangle$,
hence $\ell \times d \cdot \langle 1 \rangle \cong(d-2)\cdot\langle 1\rangle\perp\langle d/2,d/2\rangle\cong d \cdot \langle 1 \rangle$.

Conversely, if $d+1=0$ in $F$ and \(\ell\) is a similarity factor of \(d\cdot\langle 1\rangle\), then by Lemma~\ref{lem:fin_qd},
\[ \ell\times q_{d-1}\perp\langle\ell d/2\rangle\cong\ell\times((d-2)\cdot\langle1\rangle\perp\langle d/2,d/2\rangle)\cong \ell\times d\cdot\langle 1\rangle\cong d\cdot\langle 1\rangle,\]
so \(\ell\times q_{d-1}\) is a subspace of \(d\cdot\langle 1\rangle\).
Since $d\neq 0$ in $F$, $\ell\times q_{d-1}$ contains an affinely independent \(\ell\)-equilateral set of size \(d\).
By Lemma~\ref{lem:rado}, this set together with its centroid is an affinely dependent \(\ell\)-equilateral set of size \(d+1\).

\smallskip
Finally, consider statement~\ref{item:aff_dep}.\ when $d$ is odd.
Note that if $d+1=0$ in $F$, then $\ell \times q_{d-1} \cong (d-2) \cdot \langle \ell \rangle \perp \langle \ell d \rangle$ by Lemma~\ref{lem:fin_qd}, and for any $\ell\in F^*$,
\[\langle 1,1,1\rangle\cong\langle-1,-1,1\rangle\cong\langle-1,-\ell,\ell\rangle=\langle d,d\ell,\ell\rangle,\]
since $\langle \ell, -\ell \rangle \cong \langle 1, -1 \rangle$ (hyperbolic plane), and $\langle -1, -1 \rangle \cong \langle 1, 1 \rangle$ ($-1$ is in the prime subfield).

If there is an affinely dependent $\ell$-equilateral set of size $d+1$, then as in the even-dimensional case above, $d + 1 = 0$ in~$F$, and a translate of the equilateral set spans a non-degenerate hyperplane isometric to $\ell \times q_{d-1}$.
Thus,
\begin{align*}
d\cdot\langle1\rangle &\cong\ell\times q_{d-1}\perp(\ell\times q_{d-1})^\perp\\
&\cong\ell\times q_{d-1}\perp\langle d\rangle \quad\text{(comparing determinants)}\\
&\cong (d-2)\cdot\langle\ell\rangle\perp\langle\ell d,d\rangle,
\end{align*}
hence
\begin{align*}
\ell\times(d+1)\cdot\langle 1\rangle &\cong (d-2)\cdot\langle 1\rangle\perp\langle d,d\ell,\ell\rangle\cong (d+1)\cdot\langle1\rangle.
\end{align*}

Conversely, if $d+1=0$ in $F$ and $\ell$ is a similarity factor of $(d+1)\cdot\langle 1\rangle$, then
\[(d+1)\cdot\langle 1\rangle \cong (d-2)\cdot\langle1\rangle\perp\langle d,d\ell,\ell\rangle,\]
consequently
\[(d+1)\cdot\langle 1\rangle \cong \ell \times (d+1) \cdot \langle 1 \rangle \cong (d-2) \cdot \langle \ell \rangle \perp \langle \ell d, d, 1 \rangle\cong\ell\times q_{d-1}\perp\langle d,1\rangle,\]
and by Witt cancellation, $d\cdot\langle 1\rangle\cong \ell\times q_{d-1}\perp\langle d\rangle$.
This shows that $(F^d,\sum x_i^2)$ contains a hyperplane with an affinely independent $\ell$-equilateral set of size $d$.
Since $d+1=0$, Lemma~\ref{lem:rado} gives that this set, together with its centroid, is an affinely dependent $\ell$-equilateral set of size $d+1$.
\end{proof}

\begin{proof}[Proof of Theorem~\ref{cor:schoenberg_p}]
The result follows from Theorem~\ref{thm:schoenbergp}.
When $d$ is odd and $d + 1 \neq 0$ in~$F$, we take $\ell=(d+1)/2$.
When $d$ is odd and $d + 1 = 0$ in~$F$, we take $\ell=1$.
When $d$ is even, then, independent of whether $d+1=0$ or not, we also take $\ell=1$.
\end{proof}

\begin{proof}[Proof of Theorem~\ref{cor:d+2}]
Immediate from Theorems~\ref{cor:schoenberg_p} and \ref{cor:rado}.
\end{proof}

\section{Epilogue: \texorpdfstring{$\ell$}{l}-equilateral sets of size \texorpdfstring{$d$}{d}}\label{sec:size_d}
Generalizing Schoenberg's theorem to equilateral sets of size $d$ of unspecified length is trivial, as the standard basis is a $2$-equilateral set in the quadratic space $(F^d,\sum x_i^2)$.
However, we can try to generalize Theorems~\ref{thm:schoenberg0} and \ref{thm:schoenbergp} to characterize for which fields $F$ and $\ell\in F^*$ there exists an $\ell$-equilateral set of size $d$ or even of size $d-1$.
This has been done for $\ell=1$ in \cites{ElsholtzKlotz,LZ2023} and for general $\ell$ there are indirect results in \cite{LZ2023}.
For $F=\Q$, Kominers~\cite{Kominers2026}*{Theorem~1} gives a complete characterization for every positive rational~$\ell$, including sets of size~$d$ and~$d-1$.
For general fields and prescribed~$\ell$, the only case where we could find a satisfying answer is when $d$ is even.
\begin{proposition}\label{prop:exist_d_even}
Let \(d\) be even, \(F\) a field of characteristic \(\neq 2\), and
\(\ell\in F^*\). There exists an \(\ell\)-equilateral set of size \(d\)
in \((F^d,\sum_{i=1}^d x_i^2)\) if and only if \(\ell\) is a similarity factor of
\(d\cdot\langle 1\rangle\).
\end{proposition}
\begin{proof}
If $d=2$, note that there exists an $\ell$-equilateral set of size $2$ if and only if $\ell\in D(\langle 1,1\rangle)$, which holds if and only if $\langle\ell,x\rangle\cong\langle 1,1\rangle$ for some $x\in F^*$ by the representation criterion. Comparing determinants, we find that $x=\ell$ up to a square, so an equivalent condition is that $\ell\times\langle 1,1\rangle\cong\langle 1,1\rangle$.

For the rest of the proof, we assume $d\geq 4$.
We first consider the case when \(F\) has characteristic~\(0\).
By Lemma~\ref{lem:gram}, there exists an \(\ell\)-equilateral set of size \(d\) if and only if
\(\ell\times q_{d-1}\perp\langle x\rangle\cong d\cdot\langle 1\rangle\)
for some \(x\in F^*\).
By Corollary~\ref{cor:three_copies},
\(\ell\times q_{d-1}\cong(d-4)\cdot\langle\ell\rangle\perp\langle 2d\ell, 2d\ell, 2d\ell\rangle\).
Comparing determinants, we see that necessarily \(x=2d\ell\) up to a
square. Thus
\(\ell\times q_{d-1}\perp\langle 2d\ell\rangle\cong (d-4)\cdot\langle\ell\rangle\perp 4\cdot\langle 2d\ell\rangle\cong d\cdot\langle \ell\rangle\),
since \(\langle 2d,2d,2d,2d\rangle\cong\langle 1,1,1,1\rangle\).
Finally, we obtain that
\(\ell\times q_{d-1}\perp\langle 2d\ell\rangle\cong d\cdot\langle 1\rangle\)
iff \(d\cdot\langle\ell\rangle\cong d\cdot\langle 1\rangle\). 

Next we assume \(F\) has finite characteristic \(\neq 2\).
If $d\neq 0$ in $F$, again by Lemma~\ref{lem:gram} there exists an affinely independent \(\ell\)-equilateral set of size \(d\) in \(d\cdot\langle 1\rangle\) if and only if \(\ell\times q_{d-1}\perp\langle x\rangle\cong d\cdot\langle 1\rangle\) for some \(x\in F^*\), which is equivalent to
\[\ell\times d\cdot\langle 1\rangle \cong q_{d-1}\perp\langle \ell x\rangle
\cong (d-2)\cdot\langle 1\rangle\perp\langle d/2,\ell x\rangle\]
for some $x\in F^*$.
By comparing determinants, we find that necessarily $x=d/(2\ell)$ up to a square.
We conclude that there exists an affinely independent \(\ell\)-equilateral set of size \(d\) in \(d\cdot\langle 1\rangle\) if and only if $\ell\times d\cdot\langle 1\rangle \cong(d-2)\cdot\langle1\rangle\perp\langle d/2,d/2\rangle\cong d\cdot\langle 1\rangle$.

If $d=0$ in $F$, then any $\ell$-equilateral set of size $d$ contains an affinely independent subset of size $d-1$, since $d-1\neq 0$ in $F$ and Lemma~\ref{lem:rado} excludes affine dependence for such a subset.
Conversely, by the same lemma, any such subset extends by its centroid to an affinely dependent $\ell$-equilateral set of size $d$.
Thus an $\ell$-equilateral set of size $d$ exists if and only if \(\ell\times q_{d-2}\perp\langle x,y\rangle\cong d\cdot\langle1\rangle\) for some \(x,y\in F^*\), which is equivalent to \(\ell\times d\cdot\langle 1\rangle\) being isometric to
\begin{align*}
q_{d-2}\perp\langle \ell x,\ell y\rangle &\cong (d-3)\cdot\langle1\rangle\perp\langle d-1,\ell x,\ell y\rangle\\
&\cong (d-3)\cdot\langle1\rangle\perp\langle -1,\ell x,-\ell x\rangle\quad\text{(comparing determinants)}\\
&\cong (d-3)\cdot\langle1\rangle\perp\langle -1,1,-1\rangle\cong d\cdot\langle1\rangle. \qedhere
\end{align*}
\end{proof}

When $d$ is odd, we can show the following more complicated equivalences.
We omit the proofs as they are similar to that of Proposition~\ref{prop:exist_d_even}.
\begin{proposition}\label{prop:exist_d_odd}
Let $d \geq 3$ be odd, $F$ a field of characteristic $\neq 2$, and $\ell \in F^*$.
There exists an $\ell$-equilateral set of size~$d$ in $(F^d, \sum_{i=1}^d x_i^2)$ if and only if:
\begin{enumerate}
\item
$(d-1) \cdot \langle \ell \rangle \perp \langle 2\ell, d \rangle
  \cong d \cdot \langle 1 \rangle \perp \langle 2d\ell \rangle$
if $F$ has characteristic $0$,
\item
$(d-2) \cdot \langle \ell \rangle \perp \langle d\ell,\, d \rangle
  \cong d \cdot \langle 1 \rangle$
if $F$ has characteristic $p$, where $p$ is an odd prime and $p\notdivides d$,
\item
$(d-3) \cdot \langle \ell \rangle \perp \langle -2\ell,\, x,\, -2\ell x \rangle
  \cong d \cdot \langle 1 \rangle$ for some $x \in F^*$ 
if $F$ has characteristic $p$, where $p$ is an odd prime and $p\divides d$.
\end{enumerate}
\end{proposition}

\begin{problem}\label{problem}
Let $F$ be a field of characteristic not equal to $2$ and $\ell\in F^*$.
Find a simple characterization for the existence of an $\ell$-equilateral set of size $d$ in $(F^d,\sum x_i^2)$ when the dimension $d$ is odd.
\end{problem}

As mentioned in Section~\ref{section:related_work}, Kominers \cite{Kominers2026}*{Theorem~1} settled this question over~$\Q$.

\subsection*{Acknowledgements}
The author thanks the referee for suggesting a simplified proof of Lemma~\ref{lem:decomp}.

\subsection*{AI Declaration}
GPT-6 Astra was used for final proof-reading only.
The results of this paper were obtained without the use of AI.

\end{document}